\documentclass{amsart}

\usepackage{amsmath}
\usepackage{amsthm}
\usepackage{amsfonts}
\usepackage{amssymb}
\usepackage{latexsym}
\usepackage{amscd}
\usepackage{stmaryrd}

\theoremstyle{plain}
\newtheorem{thm}{Theorem}[section]
\newtheorem{cor}[thm]{Corollary}
\newtheorem{lem}[thm]{Lemma}
\newtheorem{prop}[thm]{Proposition}

\theoremstyle{plain}
\newtheorem*{conj}{Kaneko-Ohno Conjecture}

\def\proof{\noindent {\bf Proof.\;}}

\def\weight{\operatorname{wt}}
\def\depth{\operatorname{dep}}
\def\height{\operatorname{ht}}

\numberwithin{equation}{section}

\begin{document}

\title{On a conjecture of Kaneko and Ohno}

\date{}

\author{Zhong-hua Li}

\maketitle

\begin{center}
{\small Department of Mathematics, Tongji University, No. 1239 Siping Road},\\
{\small Shanghai 200092, China}\\
{\small Graduate School of Mathematical Sciences, The University of Tokyo}, \\
{\small 3-8-1 Komaba, Meguro, Tokyo 153-8914, Japan}\\
{\small E-mail address: lizhmath@gmail.com}
\end{center}

\vskip10pt

{\footnotesize
\begin{quote}

\noindent {\bf Abstract.} Let $X_0^{\star}(k,n,s)$ denote the sum of all multiple zeta-star values of weight $k$, depth $n$ and height $s$. Kaneko and Ohno conjecture that for any positive integers $m,n,s$ with $m,n\geqslant s$, the difference $(-1)^mX_0^{\star}(m+n+1,n+1,s)-(-1)^nX_0^{\star}(m+n+1,m+1,s)$ can be expressed as a polynomial of zeta values with rational coefficients. We give a proof of this conjecture in this paper.


\noindent{\bf 2010MSC}: 11M32, 33C20
\end{quote}
}


\section{Introduction}

Let $\mathbf{k}=(k_1,\ldots,k_n)$ be a sequence of positive integers with $k_1>1$,  the weight $\weight(\mathbf{k})$, depth $\depth(\mathbf{k})$ and height $\height(\mathbf{k})$ are defined by
$$\weight(\mathbf{k})=k_1+\cdots+k_n,\;\depth(\mathbf{k})=n,\;\height(\mathbf{k})=\sharp \{i\mid k_i\geqslant 2\},$$
respectively. For such a sequence $\mathbf{k}$, there are two well studied real numbers: multiple zeta value $\zeta(\mathbf{k})$ defined by
$$\zeta(\mathbf{k})=\zeta(k_1,\ldots,k_n)=\sum\limits_{m_1>\cdots>m_n>0}\frac{1}{m_1^{k_1}\cdots m_n^{k_n}},$$
and multiple zeta-star value $\zeta^{\star}(\mathbf{k})$ defined by
$$\zeta^{\star}(\mathbf{k})=\zeta^{\star}(k_1,\ldots,k_n)=\sum\limits_{m_1\geqslant \cdots\geqslant m_n\geqslant 1}\frac{1}{m_1^{k_1}\cdots m_n^{k_n}}.$$
We call the values $\zeta(\mathbf{k})$ and $\zeta^{\star}(\mathbf{k})$ with weight $\weight(\mathbf{k})$, depth $\depth(\mathbf{k})$ and height $\height(\mathbf{k})$.

The well-known Ohno-Zagier relation (\cite{ohno-zagier}) is a class of relations about the sums of multiple zeta values of fixed weight, depth and height. For integers $k,n,s$ with $k\geqslant n+s$ and $n\geqslant s\geqslant 1$, we denote by $X_0(k,n,s)$ the sum of all multiple zeta values of weight $k$, depth $n$ and height $s$. The Ohno-Zagier relation says that
$$X_0(k,n,s)\in \mathbb{Q}[\zeta(2),\zeta(3),\zeta(5),\ldots].$$
More explicitly, Ohno and Zagier gave the generating function expression
\begin{align*}
&\sum\limits_{k\geqslant n+s,n\geqslant s\geqslant 1}X_0(k,n,s)u^{k-n-s}v^{n-s}t^{s-1}\\
=&\frac{1}{uv-t}\left\{1-\exp\left(\sum\limits_{n=2}^\infty \frac{\zeta(n)}{n}(u^n+v^n-\alpha^n-\beta^n)\right)\right\},
\end{align*}
where $\alpha$ and $\beta$ are determined by $\alpha+\beta=u+v$ and $\alpha\beta=t$.
In \cite{li2010}, we showed that the Ohno-Zagier relation can be deduced from the regularized double shuffle relation. In \cite{li2008}, we generalized the concept height to $i$-height, studied sums of multiple zeta values of fixed weight, depth and general height, and expressed a kind of generating function of these sums in terms of generalized hypergeometric functions.

Similarly, we denote by $X^{\star}_0(k,n,s)$ the sum of all multiple zeta-star values of weight $k$, depth $n$ and height $s$ for integers $k,n,s$ with $k\geqslant n+s$ and $n\geqslant s\geqslant 1$.
The authors of \cite{aoki-kombu-ohno} considered a generating function $\Phi_0^{\star}(u,v,t)$ of sums $X_0^{\star}(k,n,s)$, where
$$\Phi_0^{\star}(u,v,t)=\sum\limits_{k\geqslant n+s,n\geqslant s\geqslant 1}X_0^{\star}(k,n,s)u^{k-n-s}v^{n-s}t^{2s-2}.$$
It was proved in \cite{aoki-kombu-ohno} that $\Phi_0^{\star}(u,v,t)$ can be expressed by a special value of the generalized hypergeometric function $\,_3F_2$ as
\begin{align}
\Phi_0^{\star}(u,v,t)=\frac{1}{(1-v)(1-\beta)}\,_3F_2\left({1-\beta,1-\beta+u,1\atop 2-v,2-\beta};1\right),
\label{Eq:3f2}
\end{align}
where $\alpha,\beta$ are determined by $\alpha+\beta=u+v, \alpha\beta=uv-t^2$, and the generalized hypergeometric function $\,_3F_2$ is defined as (see \cite{bailey})
$$\,_3F_2\left({\alpha_1,\alpha_2,\alpha_3\atop \beta_1,\beta_2};z\right)=\sum\limits_{n=0}^\infty
\frac{(\alpha_1)_n(\alpha_2)_n(\alpha_3)_n}{n!(\beta_1)_n(\beta_2)_n}z^n,$$
with the Pochhammer symbol $(a)_n$ given by
$$(a)_n=\frac{\Gamma(a+n)}{\Gamma(a)}=\left\{\begin{array}{ll}
1, & \text{if\;} n=0,\\
a(a+1)\cdots (a+n-1), & \text{if\;} n>0.
\end{array}\right.$$
Similarly to \cite{li2008}, the authors of \cite{aoki-ohno-wakabayashi} considered a kind of generating function of sums of multiple zeta-star values of fixed weight, depth and general height, and represented this generating function via generalized hypergeometric functions.

Since the generating function $\Phi_0^{\star}(u,v,t)$ is represented by $\,_3F_2$ as in \eqref{Eq:3f2}, it is expected that in general $X_0^{\star}(k,n,s)$ can't be written as a polynomial of zeta values with rational coefficients.  While in \cite{kaneko-ohno} Kaneko and Ohno considered some kind of duality of multiple zeta-star values, and proposed the following conjecture.
\begin{conj}[\cite{kaneko-ohno}]
For any positive integers $m,n,s$ with $m,n\geqslant s$, we have
$$(-1)^mX_0^{\star}(m+n+1,n+1,s)-(-1)^nX_0^{\star}(m+n+1,m+1,s)\in\mathbb{Q}[\zeta(2),\zeta(3),\zeta(5),\ldots].$$
\end{conj}
It was proved in \cite{kaneko-ohno} that the conjecture is true for $s=1$. Using the result of \cite{aoki-kombu-ohno} about the generating function $\Phi_0^{\star}(u,v,0)$, Yamazaki gave another proof of this case in \cite{yamazaki}. Note that the Kaneko-Ohno theorem for their conjecture in the case $s=1$ can be restated as
\begin{align}
&u\Phi_0^{\star}(-u,v,0)-v\Phi_0^{\star}(-v,u,0)
\label{Eq:height-one}\\
=&\frac{1}{u}-\frac{1}{v}+\frac{\Gamma(u+v)}{\Gamma(u)\Gamma(v)}((\Gamma(v)\Gamma(1-v))^2-(\Gamma(u)\Gamma(1-u))^2).\nonumber
\end{align}

The purpose of this paper is to give a proof of Kaneko-Ohno Conjecture. In fact, similarly to \eqref{Eq:height-one}, we give an expression of $u\Phi_0^{\star}(-u,v,t)-v\Phi_0^{\star}(-v,u,t)$ by gamma functions in Theorem \ref{Thm:main theorem}.
Our proof is based on the expression of $\Phi_0^{\star}(u,v,t)$ given in \cite{aoki-kombu-ohno}, and hence is similar to the one of Yamazaki given in \cite{yamazaki} for the special case $s=1$.

In Section 2, we state our main result and give some corollaries. In Section 3, we prepare a result about generalized hypergeometric series $\,_3F_2$. In the last section, we give the proof of the main theorem.


\section{Statement of the main result}

\subsection{Main theorem}

As in Section 1, we denote by $X_0^{\star}(k,n,s)$ the sum of all multiple zeta-star values of weight $k$, depth $n$ and height $s$ for integers $k,n,s$ with $k\geqslant n+s$ and $n\geqslant s\geqslant 1$. Let $\Phi_0^{\star}(u,v,t)$ be the generating function defined by
$$\Phi_0^{\star}(u,v,t)=\sum\limits_{k\geqslant n+s,n\geqslant s\geqslant 1}X_0^{\star}(k,n,s)u^{k-n-s}v^{n-s}t^{2s-2}.$$

For variables $u,v,t$, we define $a$ and $b$ by the conditions $a+b=-u+v$ and $ab=-uv-t^2$. Equivalently, we have
$$a,b=\frac{-u+v\pm \sqrt{(u+v)^2+4t^2}}{2}.$$
After that we define the function $A(u,v,a,b)$ by
\begin{align}
A(u,v,a,b)=\frac{1}{2\pi}\left\{\frac{\cos \pi u}{\sin \pi v}-\frac{\cos \pi v}{\sin \pi u}+\cos\pi(a-b)(\cot\pi u-\cot\pi v)\right\}.
\label{Eq:def-A}
\end{align}
Note that $A(u,v,a,b)=A(u,v,b,a)$, which shall play an important role in the proof of our main theorem.
We can express $A(u,v,a,b)$ by gamma functions as in the following lemma.

\begin{lem}
We have
\begin{align}
A(u,v,a,b)=\frac{1}{\Gamma(u+a)\Gamma(1-u-a)}\left(\frac{\Gamma(v)\Gamma(1-v)}{\Gamma(a)\Gamma(1-a)}+\frac{\Gamma(u)\Gamma(1-u)}{\Gamma(b)\Gamma(1-b)}\right),
\label{Eq:equi-A-a}
\end{align}
and
\begin{align}
A(u,v,a,b)=\frac{1}{\Gamma(u+b)\Gamma(1-u-b)}\left(\frac{\Gamma(v)\Gamma(1-v)}{\Gamma(b)\Gamma(1-b)}+\frac{\Gamma(u)\Gamma(1-u)}{\Gamma(a)\Gamma(1-a)}\right).
\label{Eq:equi-A-b}
\end{align}
\end{lem}

\proof Equation \eqref{Eq:equi-A-b} follows from equation \eqref{Eq:equi-A-a} and the fact $A(u,v,a,b)=A(u,v,b,a)$. Using the well-known reflection formula
$$\Gamma(s)\Gamma(1-s)=\frac{\pi}{\sin \pi s},$$
we find that the right-hand side of equation \eqref{Eq:equi-A-a} becomes
$$\frac{\sin \pi(u+a)}{\pi}\left(\frac{\sin\pi a}{\sin\pi v}+\frac{\sin\pi b}{\sin\pi u}\right),$$
which is equal to
$$\frac{1}{2\pi}\left(\frac{\cos\pi u-\cos\pi(v+a-b)}{\sin\pi v}+\frac{\cos\pi (u+a-b)-\cos\pi v}{\sin\pi u}\right).$$
Now it is easy to finish the proof.
\qed

The main theorem of this paper is the following theorem.

\begin{thm}\label{Thm:main theorem}
We have
\begin{align}
&u\Phi_0^{\star}(-u,v,t)-v\Phi_0^{\star}(-v,u,t)
\label{Eq:mainThm}\\
=&\frac{u-v}{ab}+A(u,v,a,b)\frac{\Gamma(a)\Gamma(1-a)\Gamma(b)\Gamma(1-b)\Gamma(u+a)\Gamma(u+b)}{\Gamma(u)\Gamma(v)}.\nonumber
\end{align}
\end{thm}

\subsection{Some remarks}

By the definition of the generating function $\Phi_0^{\star}(u,v,t)$, it is easy to see that
\begin{align}
&u\Phi_0^{\star}(-u,v,t)-v\Phi_0^{\star}(-v,u,t)
\label{Eq:difference}\\
=&\sum\limits_{m,n\geqslant s\atop s\geqslant 1}(-1)^s\left((-1)^mX_0^{\star}(m+n+1,n+1,s)-(-1)^nX_0^{\star}(m+n+1,m+1,s)\right)
\nonumber\\
&\times u^{m+1-s}v^{n+1-s}t^{2s-2}+\sum\limits_{n\geqslant s\geqslant 1}(-1)^{n+s}X_0^{\star}(n+s,s,s)(u^{n+1-s}-v^{n+1-s})t^{2s-2}.\nonumber
\end{align}
Since we have the expansion
$$\Gamma(1-x)=\exp\left(\gamma x+\sum\limits_{n=2}^\infty \frac{\zeta(n)}{n}x^n\right),$$
where $\gamma$ is the Euler's constant,
we know that Theorem \ref{Thm:main theorem} indeed implies Kaneko-Ohno Conjecture.

\begin{cor}
For any positive integers $m,n,s$ with $m,n\geqslant s$, the difference
$$(-1)^mX_0^{\star}(m+n+1,n+1,s)-(-1)^nX_0^{\star}(m+n+1,m+1,s)$$
can be expressed as a polynomial of zeta values with rational coefficients.
\end{cor}

Pay attention to the second term of the right-hand side of equation \eqref{Eq:difference}, we have another corollary.

\begin{cor}
For any positive integers $k,s$ with $k\geqslant 2s$, the sum $X_0^\star(k,s,s)$ can be expressed as a polynomial of zeta values with rational coefficients.
\end{cor}

Note that the above corollary is an immediate consequence of the symmetric sum formula for multiple zeta-star values (see \cite[Theorem 2.1]{hoffman}).

Let $t=0$ in Theorem \ref{Thm:main theorem}, we can get equation \eqref{Eq:height-one}. In fact, in this case, we can assume that $a=-u$ and $b=v$. For $A(u,v,a,b)$, we use the equivalent equation \eqref{Eq:equi-A-a}. Then using  Theorem \ref{Thm:main theorem}, we really get equation \eqref{Eq:height-one}.


\section{A result about generalized hypergeometric series $\,_3F_2$}

To prove the main theorem of this paper, we introduce the following result.

\begin{prop}\label{Prop:keyProp}
Let $a,b,c\in \mathbb{C}$ with their real parts sufficient small. We have
\begin{align}
&\,_3F_2\left({a,b,c\atop a+b, 1+c}; 1\right)
\label{Eq:keyEq}\\
=&\frac{\Gamma(a+b)\Gamma(1+c)\Gamma(1+c-a-b)}{\Gamma(a)\Gamma(b)\Gamma(1+c-a)\Gamma(1+c-b)}(\psi(1+c-b)-\psi(a)-\psi(b)-\gamma)\nonumber\\
&-\frac{\Gamma(a+b)\Gamma(1+c)\Gamma(1+c-a-b)}{\Gamma(a)\Gamma(b)\Gamma(1+c-a)\Gamma(1+c-b)}\sum\limits_{n=1}^\infty \frac{(a)_n(1-b)_n}{n n!(1+c-b)_n},
\nonumber
\end{align}
where $\psi(x)=\Gamma'(x)/\Gamma(x)$ is the digamma function.
\end{prop}

To save space, from now on we will denote the special value $\,_3F_2\left({\alpha_1,\alpha_2,\alpha_3\atop\beta_1,\beta_2};1\right)$ by  $\,_3F_2\left({\alpha_1,\alpha_2,\alpha_3\atop\beta_1,\beta_2}\right)$.

To prove the above proposition, we need two transformation formulas. The first one is (see \cite[Sec. 3.8, Eq. (1), p. 21]{bailey})
\begin{align}
\,_3F_2\left({\alpha_1,\alpha_2,\alpha_3\atop \beta_1,\beta_2}\right)
=&\frac{\Gamma(\beta_1)\Gamma(\beta_1-\alpha_1-\alpha_2)}{\Gamma(\beta_1-\alpha_1)\Gamma(\beta_1-\alpha_2)}\,_3F_2\left({\alpha_1,\alpha_2,\beta_2-\alpha_3\atop \alpha_1+\alpha_2-\beta_1+1,\beta_2}\right)
\label{Eq:trans-to-two}\\
&+\frac{\Gamma(\beta_1)\Gamma(\beta_2)\Gamma(\alpha_1+\alpha_2-\beta_1)
\Gamma(\beta_1+\beta_2-\alpha_1-\alpha_2-\alpha_3)}{\Gamma(\alpha_1)\Gamma(\alpha_2)\Gamma(\beta_2-\alpha_3)
\Gamma(\beta_1+\beta_2-\alpha_1-\alpha_2)}
\nonumber\\
&\times \,_3F_2\left({\beta_1-\alpha_1,\beta_1-\alpha_2,\beta_1+\beta_2-\alpha_1-\alpha_2-\alpha_3\atop \beta_1-\alpha_1-\alpha_2+1,\beta_1+\beta_2-\alpha_1-\alpha_2}\right),
\nonumber
\end{align}
provided that $\Re(\beta_1+\beta_2-\alpha_1-\alpha_2-\alpha_3)>0$ and $\Re(\alpha_3-\beta_1+1)>0$. The second one is (see \cite[Ex. 7, p. 98]{bailey})
\begin{align}
\,_3F_2\left({\alpha_1,\alpha_2,\alpha_3\atop \beta_1,\beta_2}\right)=\frac{\Gamma(\beta_2)\Gamma(\beta_1+\beta_2-\alpha_1-\alpha_2-\alpha_3)}{\Gamma(\beta_2-\alpha_3)\Gamma(\beta_1+\beta_2-\alpha_1-\alpha_2)}\,_3F_2
\left({\beta_1-\alpha_1,\beta_1-\alpha_2,\alpha_3\atop \beta_1,\beta_1+\beta_2-\alpha_1-\alpha_2}\right),
\label{Eq:trans-to-one}
\end{align}
provided that $\Re(\beta_1+\beta_2-\alpha_1-\alpha_2-\alpha_3)>0$ and $\Re(\beta_2-\alpha_3)>0$.

\vskip5pt

\noindent {\bf Proof of Proposition \ref{Prop:keyProp}.}
Taking a parameter $\varepsilon$, such that $|\varepsilon|$ is sufficient small, we have
$$\,_3F_2\left({a,b,c\atop a+b, 1+c}\right)=\lim\limits_{\varepsilon\rightarrow 0}\,_3F_2\left({a,b,c\atop a+b+\varepsilon, 1+c-\varepsilon}\right).$$
Now we consider the series $\,_3F_2\left({a,b,c\atop a+b+\varepsilon, 1+c-\varepsilon}\right)$.
Applying \eqref{Eq:trans-to-two}, we get
\begin{align}
\,_3F_2\left({a,b,c\atop a+b+\varepsilon, 1+c-\varepsilon}\right)=&\frac{\Gamma(a+b+\varepsilon)\Gamma(\varepsilon)}{\Gamma(a+\varepsilon)\Gamma(b+\varepsilon)}\,_3F_2\left({a,b,1-\varepsilon\atop 1-\varepsilon, 1+c-\varepsilon}\right)
\label{Eq:in-proofProp-1}\\
+&\frac{\Gamma(a+b+\varepsilon)\Gamma(1+c-\varepsilon)\Gamma(-\varepsilon)}{\Gamma(a)\Gamma(b)\Gamma(1-\varepsilon)\Gamma(1+c)}\,_3F_2\left({a+\varepsilon,
b+\varepsilon,1\atop 1+\varepsilon, 1+c}\right).
\nonumber
\end{align}
To  the first $\,_3F_2$-series in the right-hand side of \eqref{Eq:in-proofProp-1}, we apply the Gaussian summation formula (see \cite[Sec. 1.3, Eq. (1)]{bailey})
$$\sum\limits_{n=0}^\infty\frac{(\alpha_1)_n(\alpha_2)_n}{n!(\beta)_n}=
\frac{\Gamma(\beta)\Gamma(\beta-\alpha_1-\alpha_2)}{\Gamma(\beta-\alpha_1)\Gamma(\beta-\alpha_2)}$$
for $\Re(\beta-\alpha_1-\alpha_2)>0$, and apply \eqref{Eq:trans-to-one} to the second $\,_3F_2$-series in the right-hand side of \eqref{Eq:in-proofProp-1}, we obtain
\begin{align*}
&\,_3F_2\left({a,b,c\atop a+b+\varepsilon, 1+c-\varepsilon}\right)=\frac{\Gamma(1+\varepsilon)\Gamma(a+b+\varepsilon)
\Gamma(1+c-\varepsilon)\Gamma(1+c-a-b-\varepsilon)}{\varepsilon\Gamma(a+\varepsilon)\Gamma(b+\varepsilon)\Gamma(1+c-a-\varepsilon)\Gamma(1+c-b-\varepsilon)}\\
&-\frac{\Gamma(a+b+\varepsilon)\Gamma(1+c-\varepsilon)\Gamma(1+c-a-b-\varepsilon)}{\varepsilon\Gamma(a)\Gamma(b)\Gamma(1+c-a-\varepsilon)\Gamma(1+c-b)}\,_3F_2\left({\varepsilon,
1-b,a+\varepsilon\atop 1+\varepsilon, 1+c-b}\right).
\end{align*}
To the $\,_3F_2$-series in the right-hand side of the above equation, we split it into two terms as
$\sum\limits_{n=0}^\infty a_n=a_0+\sum\limits_{n=1}^\infty a_n$. Then we see that $\,_3F_2\left({a,b,c\atop a+b+\varepsilon, 1+c-\varepsilon}\right)$ is equal to
\begin{align*}
&\frac{1}{\varepsilon}\left(\frac{\Gamma(1+\varepsilon)\Gamma(a+b+\varepsilon)
\Gamma(1+c-\varepsilon)\Gamma(1+c-a-b-\varepsilon)}{\Gamma(a+\varepsilon)\Gamma(b+\varepsilon)\Gamma(1+c-a-\varepsilon)\Gamma(1+c-b-\varepsilon)}\right.\\
&-\left.\frac{\Gamma(a+b+\varepsilon)\Gamma(1+c-\varepsilon)\Gamma(1+c-a-b-\varepsilon)}{\Gamma(a)\Gamma(b)\Gamma(1+c-a-\varepsilon)\Gamma(1+c-b)}\right)\\
&-\frac{\Gamma(a+b+\varepsilon)\Gamma(1+c-\varepsilon)\Gamma(1+c-a-b-\varepsilon)}{\Gamma(a)\Gamma(b)\Gamma(1+c-a-\varepsilon)\Gamma(1+c-b)}\sum\limits_{n=1}^\infty
\frac{(a+\varepsilon)_n(1-b)_n}{(n+\varepsilon) n! (1+c-b)_n}.
\end{align*}
Finally, let $\varepsilon$ go to $0$ to finish the proof. For the first two lines of the above expression, we use L'H\^{o}pital's rule and the fact that $\psi(1)=-\gamma$.
\qed


\section{Proof of the main theorem}

In this section, we prove Theorem \ref{Thm:main theorem}.

Using the result of Aoki-Kombu-Ohno (\cite{aoki-kombu-ohno}) for the generating function $\Phi_0^{\star}(u,v,t)$, we have the following lemma.
\begin{lem}\label{Lem:ako}
Let $\alpha$ and $\beta$ be determined by $\alpha+\beta=u+v$ and $\alpha\beta=uv-t^2$. We  have
\begin{align*}
\Phi_0^{\star}(u,v,t)=&\frac{\Gamma(\beta-\alpha)\Gamma(1-\beta)\Gamma(v)\Gamma(1-v)}{\Gamma(1-\alpha)\Gamma(1+u-\alpha)\Gamma(1+\alpha-u)}\sum\limits_{n=0}^\infty\frac{(\alpha)_n(1-\beta)_n}{n!(1+\alpha-\beta)_n}\frac{\alpha-u}{n+\alpha-u}\\
&+\frac{\Gamma(\alpha-\beta)\Gamma(1-\alpha)\Gamma(v)\Gamma(1-v)}{\Gamma(1-\beta)\Gamma(1+u-\beta)\Gamma(1+\beta-u)}\sum\limits_{n=0}^\infty\frac{(\beta)_n(1-\alpha)_n}{n!(1+\beta-\alpha)_n}\frac{\beta-u}{n+\beta-u}.
\end{align*}
\end{lem}

\proof The result of Aoki-Kombu-Ohno in \cite{aoki-kombu-ohno} gives that
\begin{align*}
\Phi_0^{\star}(u,v,t)=&\frac{\Gamma(\beta-\alpha)\Gamma(1-v)}{\Gamma(1-\alpha)\Gamma(1+u-\alpha)}\int_0^1 s^{-\beta}(1-s)^{v-1}\,_2F_1\left({\alpha,\alpha-u\atop 1+\alpha-\beta};s\right)ds\\
&+\frac{\Gamma(\alpha-\beta)\Gamma(1-v)}{\Gamma(1-\beta)\Gamma(1+u-\beta)}\int_0^1 s^{-\alpha}(1-s)^{v-1}\,_2F_1\left({\beta,\beta-u\atop 1+\beta-\alpha};s\right)ds.
\end{align*}
Here $\,_2F_1\left({a,b\atop c};s\right)$ is the Gaussian hypergeometric function given by
$$\,_2F_1\left({a,b\atop c};s\right)=\sum\limits_{n=0}^\infty \frac{(a)_n(b)_n}{n!(c)_n}s^n.$$
Hence we have
\begin{align*}
&\int_0^1 s^{-\beta}(1-s)^{v-1}\,_2F_1\left({\alpha,\alpha-u\atop 1+\alpha-\beta};s\right)ds\\
=&\sum\limits_{n=0}^\infty\frac{(\alpha)_n(\alpha-u)_n}{n!(1+\alpha-\beta)_n}\int_0^1 s^{n-\beta}(1-s)^{v-1}ds\\
=&\sum\limits_{n=0}^\infty\frac{(\alpha)_n(\alpha-u)_n}{n!(1+\alpha-\beta)_n}\frac{\Gamma(1+n-\beta)\Gamma(v)}{\Gamma(1+n+v-\beta)}.
\end{align*}
Now it is easy to finish the proof.
\qed

Recall that we have defined $a$ and $b$ by
$$a+b=-u+v,\;\;ab=-uv-t^2.$$
Using the above lemma, we immediately get the following result.

\begin{lem}\label{Lem:F}
We have
$$u\Phi_0^{\star}(-u,v,t)-v\Phi_0^{\star}(-v,u,t)=F(u,v,a,b)+F(u,v,b,a),$$
where $F(u,v,a,b)$ is defined by
\begin{align*}
&\frac{\Gamma(b-a)}{\Gamma(1-u-a)\Gamma(1+u+a)}\left(\frac{u\Gamma(v)\Gamma(1-v)\Gamma(1-b)}{\Gamma(1-a)}\sum\limits_{n=0}^\infty \frac{(a)_n(1-b)_n}{n!(1+a-b)_n}\frac{u+a}{n+u+a}\right.\\
&\left.-\frac{v\Gamma(u)\Gamma(1-u)\Gamma(1+a)}{\Gamma(1+b)}\sum\limits_{n=0}^\infty \frac{(1+a)_n(-b)_n}{n!(1+a-b)_n}\frac{u+a}{n+u+a}\right).
\end{align*}
\end{lem}

Since we have
$$\sum\limits_{n=0}^\infty \frac{(a)_n(1-b)_n}{n!(1+a-b)_n}\frac{u+a}{n+u+a}=\sum\limits_{n=0}^\infty \frac{(a)_n(-b)_n}{n!(a-b)_n}\frac{(a-b)(u+a)(n-b)}{-b(n+a-b)(n+u+a)},$$
and
$$\frac{n-b}{(n+a-b)(n+u+a)}=\frac{-a}{u+b}\frac{1}{n+a-b}+\frac{v}{u+b}\frac{1}{n+u+a},$$
we get
\begin{align*}
&\sum\limits_{n=0}^\infty \frac{(a)_n(1-b)_n}{n!(1+a-b)_n}\frac{u+a}{n+u+a}\\
=&\frac{a(u+a)\Gamma(1+a-b)}{b(u+b)\Gamma(1+a)\Gamma(1-b)}
-\frac{v(a-b)}{b(u+b)}\,_3F_2\left({a,-b,u+a\atop a-b,1+u+a}\right).
\end{align*}
In the above, we have used Gaussian summation formula for Gaussian hypergeometric function at unit argument.
Similarly, we have
\begin{align*}
&\sum\limits_{n=0}^\infty \frac{(1+a)_n(-b)_n}{n!(1+a-b)_n}\frac{u+a}{n+u+a}\\
=&\frac{b(u+a)\Gamma(1+a-b)}{a(u+b)\Gamma(1+a)\Gamma(1-b)}
+\frac{u(a-b)}{a(u+b)}\,_3F_2\left({a,-b,u+a\atop a-b,1+u+a}\right).
\end{align*}
Hence we get the following lemma.

\begin{lem}\label{Lem:F-F1-F2}
We have
$$F(u,v,a,b)=F_1(u,v,a,b)+F_2(u,v,a,b),$$
where $F_1(u,v,a,b)$ is defined by
\begin{align*}
F_1(u,v,a,b)=&\frac{(u+a)\Gamma(b-a)\Gamma(1+a-b)}{(u+b)\Gamma(1-u-a)\Gamma(1+u+a)}\left(\frac{u\Gamma(v)\Gamma(1-v)}{b\Gamma(a)\Gamma(1-a)}\right.\\
&\left.-\frac{v\Gamma(u)\Gamma(1-u)}{a\Gamma(b)\Gamma(1-b)}\right),
\end{align*}
and $F_2(u,v,a,b)$ is defined by
\begin{align*}
F_2(u,v,a,b)=&\frac{uv(a-b)\Gamma(b-a)}{(u+b)\Gamma(1-u-a)\Gamma(1+u+a)}\left(\frac{\Gamma(v)\Gamma(1-v)\Gamma(-b)}{\Gamma(1-a)}\right.\\
&\left.-\frac{\Gamma(u)\Gamma(1-u)\Gamma(a)}{\Gamma(1+b)}\right)\,_3F_2\left({a,-b,u+a\atop a-b,1+u+a}\right).
\end{align*}
\end{lem}

Now we begin to compute $F_1(u,v,a,b)+F_1(u,v,b,a)$ and $F_2(u,v,a,b)+F_2(u,v,b,a)$. For $F_1(u,v,a,b)+F_1(u,v,b,a)$, we have the following result.

\begin{lem}\label{Lem:F1Duality}
The sum $F_1(u,v,a,b)+F_1(u,v,b,a)$ equals
\begin{align*}
\frac{u-v}{ab}+\frac{(a-b)uv}{ab(u+a)(u+b)}\Gamma(b-a)\Gamma(1+a-b)A(u,v,a,b).
\end{align*}
\end{lem}

\proof Using the reflection formula for gamma function, we see  that $F_1(u,v,a,b)+F_1(u,v,b,a)$ is equal to
\begin{align*}
&\Gamma(b-a)\Gamma(1+a-b)\left\{\frac{\sin \pi (u+a)}{\pi (u+b)}\left(\frac{u\sin\pi a}{b\sin\pi v}-\frac{v\sin \pi b}{a\sin\pi u}\right)\right.\\
&-\left.
\frac{\sin \pi (u+b)}{\pi (u+a)}\left(\frac{u\sin\pi b}{a\sin\pi v}-\frac{v\sin \pi a}{b\sin\pi u}\right)\right\}.
\end{align*}
The term in the brace of the above expression is
\begin{align}
&\frac{1}{2\pi(u+b)}\left(\frac{u(\cos\pi u-\cos\pi v\cos\pi(a-b)+\sin\pi v\sin \pi(a-b))}{b\sin\pi v}\right.
\label{Eq:lem-to-F1}\\
&\left.-\frac{v(\cos\pi u\cos\pi(a-b)-\sin\pi u\sin\pi(a-b)-\cos \pi v)}{a\sin\pi u}\right)
\nonumber\\
&-\frac{1}{2\pi(u+a)}\left(\frac{u(\cos\pi u-\cos\pi v\cos\pi(b-a)+\sin\pi v\sin \pi(b-a))}{a\sin\pi v}\right.
\nonumber\\
&\left.-\frac{v(\cos\pi u\cos\pi(b-a)-\sin\pi u\sin\pi(b-a)-\cos \pi v)}{b\sin\pi u}\right).
\nonumber
\end{align}
Picking up the common factors, and noting the identities
\begin{align*}
&\frac{1}{b(u+b)}-\frac{1}{a(u+a)}=\frac{v(a-b)}{ab(u+a)(u+b)},\\
&\frac{1}{a(u+b)}-\frac{1}{b(u+a)}=\frac{u(b-a)}{ab(u+a)(u+b)},\\
&\frac{u}{b(u+b)}+\frac{v}{a(u+b)}+\frac{u}{a(u+a)}+\frac{v}{b(u+a)}=\frac{2(v-u)}{ab},
\end{align*}
we see that the expression \eqref{Eq:lem-to-F1} becomes
$$\frac{uv(a-b)}{ab(u+a)(u+b)}A(u,v,a,b)+\frac{(v-u)\sin\pi(a-b)}{ab\pi},$$
which finishes the proof.
\qed

For $F_2(u,v,a,b)+F_2(u,v,b,a)$, we apply Proposition \ref{Prop:keyProp} to get the following result.

\begin{lem}\label{Lem:F2Duality}
The sum $F_2(u,v,a,b)+F_2(u,v,b,a)$ equals
\begin{align*}
&\frac{(b-a)uv}{ab(u+a)(u+b)}\Gamma(b-a)\Gamma(1+a-b)A(u,v,a,b)\\
&+A(u,v,a,b)\frac{\Gamma(a)\Gamma(1-a)\Gamma(b)\Gamma(1-b)\Gamma(u+a)\Gamma(u+b)}{\Gamma(u)\Gamma(v)}.
\end{align*}
\end{lem}

\proof Applying Proposition \ref{Prop:keyProp} to the $\,_3F_2$-series in $F_2(u,v,a,b)$, we find that $F_2(u,v,a,b)$ becomes
\begin{align*}
&\frac{(a-b)\Gamma(a-b)\Gamma(b-a)\Gamma(u+b)}{\Gamma(u)\Gamma(v)\Gamma(1-u-a)}\left(\frac{\Gamma(v)\Gamma(1-v)}{\Gamma(a)\Gamma(1-a)}-\frac{\Gamma(u)\Gamma(1-u)}{\Gamma(-b)\Gamma(1+b)}\right)\\
&\times \left(\psi(1+v)-\psi(a)-\psi(-b)-\gamma-\sum\limits_{n=1}^\infty\frac{(a)_n(1+b)_n}{nn!(1+v)_n}\right),
\end{align*}
which is just
\begin{align*}
&A(u,v,a,b)\frac{\Gamma(b-a)\Gamma(1+a-b)\Gamma(u+a)\Gamma(u+b)}{\Gamma(u)\Gamma(v)}\\
&\times \left(\psi(1+v)-\psi(a)-\psi(-b)-\gamma-\sum\limits_{n=1}^\infty\frac{(a)_n(1+b)_n}{nn!(1+v)_n}\right).
\end{align*}
Hence using the fact $A(u,v,a,b)=A(u,v,b,a)$, we find $F_2(u,v,a,b)+F_2(u,v,b,a)$ becomes
\begin{align*}
&A(u,v,a,b)\frac{\Gamma(b-a)\Gamma(1+a-b)\Gamma(u+a)\Gamma(u+b)}{\Gamma(u)\Gamma(v)}\\
&\times \left\{\sum\limits_{n=1}^\infty\left(\frac{(1+a)_n(b)_n}{nn!(1+v)_n}-\frac{(a)_n(1+b)_n}{nn!(1+v)_n}\right)+\psi(b)-\psi(-b)+\psi(-a)-\psi(a)\right\}.
\end{align*}
It is easy to see that
$$\sum\limits_{n=1}^\infty\left(\frac{(1+a)_n(b)_n}{nn!(1+v)_n}-\frac{(a)_n(1+b)_n}{nn!(1+v)_n}\right)=\frac{b-a}{ab}\sum\limits_{n=1}^\infty \frac{(a)_n(b)_n}{n!(1+v)_n},$$
which equals
$$\frac{b-a}{ab}\frac{\Gamma(1+u)\Gamma(1+v)}{\Gamma(1+u+a)\Gamma(1+u+b)}+\frac{1}{b}-\frac{1}{a}$$
by Gaussian summation formula.
Applying the formulas
$$\psi(-x)-\psi(x)-\frac{1}{x}=\pi\cot\pi x,$$
and
$$\pi \cot\pi a-\pi\cot\pi b=\frac{\Gamma(a)\Gamma(1-a)\Gamma(b)\Gamma(1-b)}{\Gamma(b-a)\Gamma(1+a-b)},$$
we finish the proof.
\qed

\vskip5pt

\noindent {\bf Proof of Theorem \ref{Thm:main theorem}.} Theorem \ref{Thm:main theorem} follows from Lemma \ref{Lem:F}, Lemma \ref{Lem:F-F1-F2}, Lemma \ref{Lem:F1Duality} and Lemma \ref{Lem:F2Duality}.
\qed


\vskip5pt

\noindent {\bf Acknowledgements.} This work was partially supported by the National Natural
Science Foundation of China (grant no. 11001201), the Program for
Young Excellent Talents in Tongji University (grant no. 2009KJ065) and
the Japan Society for the Promotion of Science postdoctoral Fellowship for Foreign Researchers.



\begin{thebibliography}{99}

\bibitem{aoki-kombu-ohno} T. Aoki, Y. Kombu and Y. Ohno, A generating function for sums of multiple zeta values and its applications, Proc. Amer. Math. Soc. 136(2008), 387-395.

\bibitem{aoki-ohno-wakabayashi} T. Aoki, Y. Ohno and N. Wakabayashi, On generating functions of multiple zeta values and generalized hypergeometric functions, Manuscripta Math. 134(2011), 139-155.

\bibitem{bailey} W. N. Bailey, Generalized hypergeometric series, Cambridge University Press, Cambridge, 1935.

\bibitem{hoffman} M. E. Hoffman, Multiple harmonic series, Pacific J. Math. 152(2) (1992), 275-290.

\bibitem{kaneko-ohno} M. Kaneko and Y. Ohno, On a kind of duality of multiple zeta-star values, Int. J. Number Theory 6(8)(2010), 1927-1932.

\bibitem{li2008} Z. Li, Sum of multiple zeta values of fixed weight, depth and $i$-height, Math. Z. 258(1)(2008), 133-142.

\bibitem{li2010} Z. Li, Regularized double shuffle and Ohno-Zagier relations of multiple zeta values, preprint, 2010.

\bibitem{ohno-zagier} Y. Ohno and D. Zagier, Multiple zeta values of fixed weight, depth, and  height, Indag. Math. 12(4) (2001), 145-152.

\bibitem{yamazaki} C. Yamazaki, On the duality for multiple zeta-star values of height one, Kyushu J. Math. 64(1) (2010), 145-152.


\end{thebibliography}
\end{document}